\PassOptionsToPackage{hyperfootnotes=false}{hyperref}
\documentclass[a4paper,fleqn]{cas-sc}
\usepackage[utf8]{inputenc}
\usepackage[numbers,sort&compress]{natbib}
\usepackage{amsthm,mathtools,microtype}
\hypersetup{hidelinks,
 pdftitle={Domain-uniform integration lower bounds and local information}}
\ExplSyntaxOn
\keys_set:nn { stm / mktitle } { nologo = true }
\cs_new_eq:NN \PrintNonemptyORCID \printorcid
\RenewDocumentCommand \printorcid { }
  { \seq_if_empty:NF \g_stm_orcid_seq { \PrintNonemptyORCID } }
\ExplSyntaxOff
\numberwithin{equation}{section}
\newtheorem{theorem}{Theorem}[section]
\newtheorem{lemma}[theorem]{Lemma}
\newtheorem{proposition}[theorem]{Proposition}
\newtheorem{corollary}[theorem]{Corollary}
\theoremstyle{definition}

\theoremstyle{remark}
\newtheorem{remark}[theorem]{Remark}
\newcommand{\R}{\mathbb R}
\newcommand{\N}{\mathbb N}
\newcommand{\dd}{\,\mathrm d}
\newcommand{\e}{\mathrm e}
\newcommand{\Lip}{\operatorname{Lip}}
\newcommand{\dist}{\operatorname{dist}}
\newcommand{\op}{\mathrm{op}}
\newcommand{\F}{\mathcal F}
\newcommand{\Fn}{\mathcal F(1,d^{-1/2},d^{-1})}
\newcommand{\CD}{\mathcal C_D}
\newcommand{\B}{\mathcal B_b(\R^d)}
\newcommand{\vol}{\operatorname{vol}_d}
\newcommand{\one}{\mathbf 1}
\newcommand{\non}{\mathrm{non}}
\newcommand{\loc}{\mathrm{loc}}
\newcommand{\vg}{\mathrm{vg}}
\newcommand{\val}{\mathrm{val}}

\begin{document}
\shorttitle{Domain-uniform integration and local information}
\title[mode=title]{Domain-uniform integration lower bounds and local information}

\author[1]{Yutong Zhang}
\shortauthors{Yutong Zhang}
\cormark[1]
\ead{yutongzhang@stu.scu.edu.cn}

\affiliation[1]{
  organization={School of Mathematics, Sichuan University},
  city={Chengdu},
  postcode={610065},
  country={China}
}

\cortext[1]{Corresponding author}

\begin{abstract}
We study deterministic numerical integration on open subsets of
$\mathbb R^d$ of volume one. Let $e_n$ be the minimal worst-case error
using at most $n$ function values. For twice continuously differentiable
functions with the normalized Lipschitz bounds proposed by Novak, we prove
$e_n\geq v_d^{-2/d}n^{-2/d}/(4(d+2))\geq(8\pi\mathrm e)^{-1}n^{-2/d}$,
where $v_d$ is the volume of the Euclidean unit ball. This resolves
Novak's conjecture: the lower-bound constant is independent of the
dimension and the shape of the domain. A midpoint rule on the cube
has the same rate, up to an absolute constant.
The proof uses smooth functions that vanish near the sample points
and are large elsewhere. The bound therefore remains valid when each
sample supplies all local information about the function.
For scalar linear targets continuous in the uniform norm, we prove a
second result on the global smooth class with prescribed function,
gradient, and Hessian bounds: adaptive local information has the same
minimal error as values and gradients at points chosen in advance.
We give extensions to function recovery, finite integration measures,
and anisotropic derivative bounds.
\end{abstract}
\maketitle
\hypersetup{pdfauthor={Yutong Zhang},
 pdfsubject={Numerical integration and information-based complexity}}

\section{Introduction}\label{sec:intro}

How many sample points are needed to integrate a smooth function on a
domain of volume one? For a fixed domain, many smoothness classes have
known error rates. Here the domain may change with the dimension and
need not be convex, connected, or bounded. We prove a lower bound that
holds for every such open domain, with the same positive constant in
every dimension.

Let $D\subset\R^d$ be open with $|D|=1$, and write
$I_Df=\int_Df(x)\dd x$. All norms on $\R^d$ are Euclidean.
For a function $g$ on $D$, let
\[
 \Lip_Dg=\sup_{\substack{x,y\in D\\x\ne y}}
                  \frac{|g(x)-g(y)|}{|x-y|}.
\]
This definition compares all pairs of points, including points in
different connected components. For a unit vector $\theta$, write
$D^\theta f=\nabla f\cdot\theta$. Consider the class
\begin{equation}
\label{eq:novakclass}
 \widetilde C^2(D)=
 \left\{f\in C^2(D):
 \begin{aligned}
  &\|f\|_{\infty,D}\leq1,\qquad \Lip_D f\leq d^{-1/2},\\
  &\Lip_D(D^\theta f)\leq d^{-1}
                      \quad\text{for every }|\theta|=1
 \end{aligned}\right\}.
\end{equation}
Let $e_n(\widetilde C^2(D))$ denote the infimum of the worst-case
integration errors of deterministic algorithms using at most $n$
function values at points of $D$. Thus
\[
 e_n(\widetilde C^2(D))
 =\inf_{\mathcal A}\sup_{f\in\widetilde C^2(D)}
       |I_Df-\mathcal A(f)|.
\]
The sample points may depend on earlier values; this is called adaptive
sampling. The output may be any function of the observed data.

Novak asked whether an absolute constant $C>0$ satisfies
\begin{equation}
\label{eq:conjecture}
 e_n(\widetilde C^2(D))\geq Cn^{-2/d}
 \qquad(d,n\in\N,\ D\subset\R^d\text{ open},\ |D|=1);
\end{equation}
see \cite[Section 3, Example 2]{Novak2020} and
\cite[Section 4, Example 1]{Novak2024}.
Let $v_d=\vol(B_2^d)$ be the volume of the Euclidean unit ball.

\begin{theorem}\label{thm:main}
For every $d,n\in\N$ and every open $D\subset\R^d$ with $|D|=1$,
\begin{equation}
\label{eq:main}
 e_n(\widetilde C^2(D))\geq\frac{v_d^{-2/d}}{4(d+2)}n^{-2/d}
 \geq\frac1{8\pi\e}n^{-2/d}.
\end{equation}
\end{theorem}

The same lower bound holds for globally defined $C^\infty$ functions
satisfying
\begin{equation}
\label{eq:normalizedglobal}
 \|f\|_\infty\leq1,\qquad
 \|\nabla f\|_\infty\leq d^{-1/2},\qquad
 \|\nabla^2f\|_{\op,\infty}\leq d^{-1},
\end{equation}
even when an adaptive algorithm samples anywhere in $\R^d$ and receives
local information: two functions that agree on a neighborhood of the
sample point return the same data. Section~\ref{sec:models} gives the
precise information model.
On the unit cube, function values give
\begin{equation}
\label{eq:cubeintro}
 e_n\bigl(\widetilde C^2((0,1)^d)\bigr)
 \leq\frac16n^{-2/d}.
\end{equation}
Thus the order $n^{-2/d}$ cannot be improved in a bound uniform over all
volume-one domains.

Here $\|\nabla^2f\|_{\op,\infty}$ is the supremum of the operator norm
of the Hessian. The global functions in \eqref{eq:normalizedglobal}
restrict to members of \eqref{eq:novakclass}, even when $D$ is disconnected.
This information includes the function value and all derivatives at the
sample point.

The proof has two steps. First, run an algorithm on the zero function
and record its sample points. Construct a nonnegative smooth function
that is zero near these points and large away from them. The algorithm
returns the same answer for this function and its negative, so one of
the two errors is at least the integral of the constructed function.
Second, estimate that integral using the fact that $n$ balls of radius
$r$ cover at most $nv_dr^d$ units of volume. This estimate uses the
volume of $D$ and no further geometric property.

Earlier lower bounds for integration of smooth functions include
\cite{HNUW2014a,HNUW2014b}; results for general domains appear in
\cite{HPU2019}. The distance construction in \cite[Section 4.1]{HNUW2014b} uses the
convex hull of the sample points. We use the distance to the points
themselves and smooth it while keeping the first two derivative bounds
independent of the dimension. For sampling rates on fixed domains,
see \cite{NT2006}. General background on information-based complexity
is given in \cite{NW2008}.

\subsection{Further results and organization}
For $A,L,H>0$, we work with the global class
\begin{equation}
\label{eq:globalclass}
 \F=\mathcal F(A,L,H)=\left\{f\in C^\infty(\R^d):
 \|f\|_\infty\leq A,\quad
 \|\nabla f\|_\infty\leq L,\quad
 \|\nabla^2f\|_{\op,\infty}\leq H\right\}.
\end{equation}
We suppress $d,A,L,H$ in the notation when they are fixed. Norms
without a domain subscript are taken over $\R^d$.
Section~\ref{sec:models} states the smooth-function construction used in
the lower bounds. Its proof is in Appendix~\ref{app:construction} and
uses quadratic regularization followed by convolution. This
regularization goes back to Lasry and Lions \cite{LL1986} and Attouch
and Az\'e \cite{AA1993}; see also \cite{Bernard2010,ALM2018,ZCO2015}.

For arbitrary $A,L,H$, Section~\ref{sec:mainproof} gives a lower bound
by a universal constant times
\[
 \min\{A,L\sqrt d\,n^{-1/d},Hd\,n^{-2/d}\}.
\]
The three terms come from the bounds on the function, gradient, and
Hessian. Midpoint rules on a cube give the same order.

Section~\ref{sec:operator} compares local information with function
values and gradients. For integration against a positive measure and for function recovery,
an explicit function of the distance to the nearest sample point
estimates the minimal error within a factor of two. A positive-operator formulation covers both
problems. For scalar linear targets
continuous in the uniform norm, adaptive local information has the same
minimal error as values and gradients at points chosen in advance. The proof approximates functions whose
values and gradients vanish at the sample points by functions that
vanish on neighborhoods of those points.
Related recovery problems from values and gradients are studied in
\cite{BBS2010,BS2011,Borodachov2022}. Distance estimates for recovery in
other function classes appear in \cite{KS2024,KS2025}, and a survey of
optimal sampling is given in \cite{KU2026}.

Section~\ref{sec:extensions} applies these estimates to $L_q$ recovery,
finite measures, Gaussian integration, and anisotropic smoothness.
Section~\ref{sec:scope} discusses optimal constants and randomized
algorithms.

\section{Information and the test functions}\label{sec:models}

\subsection{Information and the lower-bound argument}
A sample at $p$ may return the value and gradient
$(f(p),\nabla f(p))$. We also allow a stronger kind of information:
two functions give the same \emph{complete local information} at $p$
exactly when they agree throughout some neighborhood of $p$.
The neighborhood may depend on the two functions. In particular, the
model does not specify a positive observation radius.

On the smooth class $\F$, complete local information determines every
derivative at $p$. To see
why it can contain more than those derivatives, consider the function
$\exp(-1/t^2)$ for $t\ne0$, extended by zero at $t=0$. Every derivative
at zero vanishes, although the function is positive at every other point.
For our lower bounds, the relevant fact is simpler: a function that is
zero near $p$ gives the same local information as the zero function.

Each sampled point has unit cost, regardless of which of these data is
returned. An adaptive algorithm may choose its next point and its stopping
time from the data already obtained. A nonadaptive algorithm chooses its
points in advance. In both cases, an algorithm with budget $n$ must stop
after at most $n$ samples on every input.

Let $\mathcal G$ be a function class containing zero and satisfying
$f\in\mathcal G\Rightarrow-f\in\mathcal G$.
Let $S$ be the quantity to be recovered: for integration, $S=I_D$;
for function recovery, $Sf$ is the function itself in the required norm.
More generally, $S$ may be a linear map from the span of $\mathcal G$
into a real normed space $Y$. The worst-case error of an algorithm
$\mathcal A$ is
\[
 e(\mathcal A,S,\mathcal G)
 =\sup_{f\in\mathcal G}\|Sf-\mathcal A(f)\|_Y.
\]
The corresponding minimal error is the infimum over all permitted
algorithms. Superscripts $\loc$ and $\vg$ indicate local information and
values with gradients. Adaptive sampling is allowed unless the
superscript $\non$ specifies points chosen in advance. The superscript
$\val$ indicates function values alone.
For a fixed finite set $P$, we write $e_S^{\mathcal I}(P)$ for
the minimal error using information of type $\mathcal I$ at $P$. These fixed-point errors refer to $\F$ unless another class is stated.

To obtain a lower bound, run the algorithm on zero and record its sample
points. We then construct $f$ so that the data at those points are zero.
The algorithm follows the same steps on $f$, $-f$, and zero. It therefore
returns one answer for the two targets $Sf$ and $-Sf$, whose distance is
$2\|Sf\|_Y$.

\begin{lemma}\label{lem:zero-data}
Run $\mathcal A$ on the zero function and let $P$ be the set of points
it samples. If $f,-f\in\mathcal G$ give the same data as zero at every
point of $P$, then
\[
 e(\mathcal A,S,\mathcal G)\geq\|Sf\|_Y.
\]
\end{lemma}
\begin{proof}
At each step, the data obtained from $f$, $-f$, and zero agree. Induction
therefore shows that the algorithm chooses the same points, stops at the
same step, and returns the same answer $\mathcal A(0)$ on all three
inputs. Since $S(-f)=-Sf$,
\begin{align*}
 2\|Sf\|_Y
 &=\|(Sf-\mathcal A(0))-(-Sf-\mathcal A(0))\|_Y\\
 &\leq\|Sf-\mathcal A(0)\|_Y+\|-Sf-\mathcal A(0)\|_Y
 \leq2e(\mathcal A,S,\mathcal G).
\end{align*}
\end{proof}

Throughout, $n\geq1$. Repeated sample points can be removed. If an
algorithm takes no samples on zero, we may choose an arbitrary point and
require the function used in the lower bound to vanish near that point.
Thus all infima over sample sets below may be taken over $1\leq|P|\leq n$.
When sample locations are restricted, every information model in a
comparison uses the same nonempty set of admissible points. By default,
samples may lie anywhere in $\R^d$ for $\F$, and must lie in $D$ for
functions defined only on $D$.

\subsection{A function of the distance}
Fix $A,L,H>0$. Let $K\subset\R^d$ be closed and nonempty, and write
$\rho_K(x)=\dist(x,K)$. The application is $K=P$, the finite set of
sample points. We want a function in $\F$ that is zero near $K$ and
large elsewhere.

Distance gives a natural starting point. Near $K$, a function with zero
value and gradient can grow quadratically under a Hessian bound. Farther
away, the gradient bound limits its growth to linear growth, and the
function bound stops it at height $A$. We use the following function of
one variable:
\[
 b(r)=\min\{A,q(r)\},\qquad
 q(r)=
 \begin{cases}
 Hr^2/4,&0\leq r\leq 2L/H,\\
 Lr-L^2/H,&r\geq 2L/H.
 \end{cases}
\]
For the normalized bounds $A=1$, $L=d^{-1/2}$, and $H=d^{-1}$
of Theorem~\ref{thm:main}, this reduces to
\[
 b(r)=\min\{1,r^2/(4d)\}.
\]
Thus the main integration theorem only uses a quadratic function capped
at height one.

The coefficient $H/4$ leaves room for the smoothing step below.
The function $b$ itself need not be smooth at the point where it reaches
$A$. We shall use it as a lower bound for the smooth test functions.
For comparison with Taylor's formula, consider
\[
 \min\{A,Lr,Hr^2/2\}.
\]
The function $q$ becomes linear at $2L/H$, and $b$ reaches $A$ at
\begin{equation}
\label{eq:radii}
 R=\begin{cases}
 2\sqrt{A/H},&A\leq L^2/H,\\
 A/L+L/H,&A\geq L^2/H.
 \end{cases}
\end{equation}
The two expressions agree when $A=L^2/H$.

\begin{lemma}\label{lem:profile}
The function $b$ is nonnegative, nondecreasing, bounded by $A$, and
$L$-Lipschitz. It satisfies
\begin{align}
 b(r)&\leq Hr^2/4,\label{eq:bquad}\\
 \tfrac12\min\{A,Lr,Hr^2/2\}&\leq b(r)\leq\min\{A,Lr,Hr^2/2\},\label{eq:bchi}\\
 \tfrac14\min\{A,Lr,Hr^2\}&\leq b(r),\label{eq:bpsi}\\
 b(tr)&\geq t^2b(r)\qquad(0\leq t\leq1),\label{eq:bscale}\\
 b'(r)&=\min\{Hr/2,L\}\one_{(0,R)}(r)
 \quad\text{for almost every }r>0.\label{eq:bderiv}
\end{align}
\end{lemma}
The proof is given in Appendix~\ref{app:profile}.

\subsection{Smooth functions vanishing near a closed set}\label{sec:lift}
The distance $\rho_K$ can have corners, so $b(\rho_K)$ need not be
smooth. The next theorem replaces it by a smooth function satisfying
the bounds defining $\F$. The replacement is zero on a
neighborhood of $K$, and its value falls below $b(\rho_K)$ by at most
$L(\delta+\eta)$. Thus the loss can be made arbitrarily small.

\begin{theorem}\label{thm:lift}
Let $K\subset\R^d$ be closed and nonempty, and let $0<\eta<\delta$.
There exists $f_{\delta,\eta}\in\F$ such that
\begin{align}
 0&\leq f_{\delta,\eta}\leq A,\notag\\
 f_{\delta,\eta}(x)&=0
 \quad\text{if }\rho_K(x)\leq\delta-\eta,\label{eq:liftflat}\\
 f_{\delta,\eta}(x)&\geq
 b\bigl((\rho_K(x)-\delta-\eta)_+\bigr)\notag\\
 &\geq b(\rho_K(x))-L(\delta+\eta).
 \label{eq:liftprofile}
\end{align}
Here $t_+=\max\{t,0\}$. The constants do not depend on the dimension or
on $K$.
\end{theorem}

The proof is given in Appendix~\ref{app:construction}.

\begin{remark}\label{rem:buffer}
In the lower-bound proofs, take $\eta=\delta/2$ and let
$\delta\downarrow0$. The loss in \eqref{eq:liftprofile} is then
$3L\delta/2$. For every positive $\delta$, the test function remains
zero within distance $\delta/2$ of $K$ and is therefore indistinguishable
from zero at the sample points in $K$.
\end{remark}

\section{Integration on arbitrary domains}\label{sec:mainproof}

\subsection{Functions on an open set}
We first transfer the global construction to the class on $D$.
For fixed $A,L,H>0$, write
\begin{equation*}
 \CD=\mathcal C_D(A,L,H)=\left\{f\in C^2(D):
 \begin{aligned}
 &\|f\|_{\infty,D}\leq A,\quad\Lip_D f\leq L,\\
 &\Lip_D(D^\theta f)\leq H\quad(|\theta|=1)
 \end{aligned}\right\}.
\end{equation*}
The restriction of every $f\in\F$ belongs to $\CD$.
Indeed, for any $x,y\in\R^d$,
\begin{align*}
 f(y)-f(x)
 &=\int_0^1\langle\nabla f(x+t(y-x)),y-x\rangle\dd t,
 \\
 D^\theta f(y)-D^\theta f(x)
 &=\int_0^1\theta^T\nabla^2f(x+t(y-x))(y-x)\dd t.
\end{align*}
These identities are available even if the segment from $x$ to $y$
leaves $D$, because $f$ is defined on all of $\R^d$. They give
\begin{equation}
\label{eq:restrictionLip}
 |f(y)-f(x)|\leq L|y-x|,\qquad
 |D^\theta f(y)-D^\theta f(x)|\leq H|y-x|.
\end{equation}
This proves the required Lipschitz bounds on $D$.

\begin{proposition}\label{prop:integralprofile}
Let $D$ be measurable with finite volume $V>0$, and let $\mathcal A$ be a
deterministic algorithm on $\F$ using local information at at most $n$ points.
Let $P$ be the set of points it samples on zero, adding one arbitrary
admissible point if this set is empty. Then
\begin{equation*}
 e(\mathcal A,I_D,\F)\geq\int_D b(\rho_P(x))\dd x.
\end{equation*}
\end{proposition}
\begin{proof}
Use Theorem~\ref{thm:lift} with $K=P$ and $\eta=\delta/2$.
The functions $f_{\delta,\delta/2}$ and $-f_{\delta,\delta/2}$ vanish
in a neighborhood of every point of $P$. By
Lemma~\ref{lem:zero-data} and nonnegativity,
\begin{align*}
 e(\mathcal A,I_D,\F)
 &\geq\int_D f_{\delta,\delta/2}(x)\dd x\\
 &\geq\int_D b(\rho_P(x))\dd x-\frac32LV\delta.
\end{align*}
The left side is independent of $\delta$. Letting $\delta\downarrow0$
proves the claim.
\end{proof}

\subsection{The volume estimate}
The integration estimate depends on how much of $D$ lies far from the
sample points. At most $n$ balls of radius $r$ cover a volume of at most
$nv_dr^d$.
Thus, for $|D|=1$ and $|P|\leq n$,
\begin{equation}
\label{eq:tubevolume}
 |D\cap\{\rho_P\leq r\}|
 \leq\left|\bigcup_{p\in P}\overline B(p,r)\right|
 \leq nv_dr^d.
\end{equation}
The cube $[-d^{-1/2},d^{-1/2}]^d$ lies in the closed Euclidean unit ball, hence
\begin{equation}
\label{eq:lowerball}
 v_d\geq(2/\sqrt d)^d,
 \qquad (nv_d)^{-1/d}\leq\frac{\sqrt d}{2}n^{-1/d}\leq\frac{\sqrt d}{2}.
\end{equation}
For the normalized parameters in \eqref{eq:normalizedglobal}, the profile
simplifies to
$b(r)=\min\{1,r^2/(4d)\}$, with $R=2\sqrt d$.
By \eqref{eq:lowerball}, $b(r)=r^2/(4d)$ throughout $[0,(nv_d)^{-1/d}]$.
For every $x$, monotonicity of $b$ gives
\[
 b(\rho_P(x))\geq b(\min\{\rho_P(x),(nv_d)^{-1/d}\})
 =\int_0^{(nv_d)^{-1/d}}\frac r{2d}\,
          \one_{\{\rho_P(x)>r\}}\dd r.
\]
Integrate this identity over $D$ and use Tonelli's theorem to exchange
the order of integration. Equation~\eqref{eq:tubevolume} then gives
\begin{align}
 \int_D b(\rho_P(x))\dd x
 &\geq\frac1{2d}\int_0^{(nv_d)^{-1/d}}r\,|D\cap\{\rho_P>r\}|\dd r\notag\\
 &\geq\frac1{2d}\int_0^{(nv_d)^{-1/d}}r(1-nv_dr^d)\dd r\notag\\
 &=\frac1{2d}\left(\frac{(nv_d)^{-2/d}}{2}
                 -\frac{nv_d(nv_d)^{-(d+2)/d}}{d+2}\right)\notag\\
 &=\frac{(nv_d)^{-2/d}}{4(d+2)}
 =\frac{v_d^{-2/d}}{4(d+2)}n^{-2/d}.
 \label{eq:mainintegral}
\end{align}
The estimate holds for every measurable $D$ of volume one.

\subsection{A bound independent of the dimension}
\begin{lemma}\label{lem:volume}
For every $d\in\N$,
\begin{equation}
\label{eq:upperball}
 v_d\leq\left(\frac{2\pi\e}{d+2}\right)^{d/2}.
\end{equation}
\end{lemma}
\begin{proof}
Integration in polar coordinates gives
$\frac1{v_d}\int_{B_2^d}|x|^2\dd x =d\int_0^1r^{d+1}\dd r=\frac d{d+2}$.
Since $u\mapsto\exp(-au)$ is convex for $a>0$, Jensen's inequality and
then the Gaussian integral imply
\begin{equation*}
 v_d\exp\!\left(-\frac{ad}{d+2}\right)
 \leq\int_{B_2^d}\e^{-a|x|^2}\dd x
 \leq\int_{\R^d}\e^{-a|x|^2}\dd x
 =\left(\frac\pi a\right)^{d/2}.
\end{equation*}
Choose $a=(d+2)/2$ and rearrange.
\end{proof}

\begin{proof}[Proof of the lower bound in Theorem~\ref{thm:main}]
Apply Proposition~\ref{prop:integralprofile}, \eqref{eq:mainintegral}, and
Lemma~\ref{lem:volume} to every algorithm on the globally smooth subclass:
\begin{equation*}
 e(\mathcal A,I_D,\Fn)
 \geq\frac{v_d^{-2/d}}{4(d+2)}n^{-2/d}
 \geq\frac1{8\pi\e}n^{-2/d}.
\end{equation*}
If the algorithm makes no query on the zero input, the constant
functions $1$ and $-1$ give error at least one. Taking the infimum over
algorithms proves the lower bound for local information. By
\eqref{eq:restrictionLip}, the restrictions of these inputs belong to
$\widetilde C^2(D)$. Every point-value algorithm has no more information
than an algorithm using local information, so the same bound proves
\eqref{eq:main}.
\end{proof}

The explicit constants are $1/48$ for $d=1$ and $1/(16\pi)$ for $d=2$, and
$v_d^{-2/d}/(4(d+2))\to1/(8\pi\e)$ as $d\to\infty$.
The limit follows from
$v_d=\pi^{d/2}/\Gamma(1+d/2)$ and Stirling's formula
$\Gamma(1+t)\sim\sqrt{2\pi t}(t/\e)^t$.

\subsection{The cube upper bound and the curse of dimensionality}
Partition $(0,1)^d$ into $m^d$ equal subcubes, and let $c_Q$ be their
centers. The midpoint formula is
$M_mf=m^{-d}\sum_Qf(c_Q)$.
At an interior point, difference quotients of $D^\theta f$ along any
unit vector $v$ give
$|\theta^T\nabla^2f(x)v|\leq d^{-1}$. Taking the supremum over
$\theta$ and $v$ yields $\|\nabla^2f(x)\|_\op\leq d^{-1}$.
The line segment from a cell center to a point of that cell stays in
the cube, so Taylor's theorem gives
\begin{equation*}
 |f(x)-f(c_Q)-\langle\nabla f(c_Q),x-c_Q\rangle|
 \leq\frac1{2d}|x-c_Q|^2.
\end{equation*}
The linear term integrates to zero and
$\frac1{|Q|}\int_Q|x-c_Q|^2\dd x =\frac d{12m^2}$.
Therefore
$|I_Df-M_mf|\leq\frac1{2d}\frac d{12m^2} =\frac1{24m^2}$.
For $m=\lfloor n^{1/d}\rfloor\geq n^{1/d}/2$, this proves
\eqref{eq:cubeintro}.

Let $n(\varepsilon,d,D)$ be the least query budget achieving error at
most $\varepsilon$. The initial error is one, and
\begin{equation*}
 0<\varepsilon<\frac1{8\pi\e}
 \quad\Longrightarrow\quad
 n(\varepsilon,d,D)\geq
 \left(\frac1{8\pi\e\,\varepsilon}\right)^{d/2}.
\end{equation*}
Thus the required number of sample points grows exponentially with $d$
for every such fixed $\varepsilon$, uniformly in $D$. The same conclusion
holds for local information on the global class
\eqref{eq:normalizedglobal}.

\subsection{General bounds on the function and its derivatives}\label{sec:parameters}

\begin{theorem}\label{thm:three}
Let $d,n\in\N$, $A,L,H>0$, and $D\subset\R^d$ be open with $|D|=1$.
Then
\begin{align}
 e_n(I_D,\CD)
 &\geq\frac d{d+2}b\bigl((nv_d)^{-1/d}\bigr)\label{eq:exactthree}\\
 &\geq\frac1{8\pi\e}
 \min\{A,L\sqrt d\,n^{-1/d},Hd\,n^{-2/d}\}.
 \label{eq:three}
\end{align}
The two lower bounds also hold on $\F$ for adaptive local information,
with sample points anywhere in $\R^d$. On $D=(0,1)^d$, point values give
\begin{equation}
\label{eq:threeupper}
 e_n(I_D,\CD)
 \leq\min\left\{A,\frac{L\sqrt d}{\sqrt3}n^{-1/d},
                       \frac{Hd}{6}n^{-2/d}\right\}.
\end{equation}
\end{theorem}
\begin{proof}
Since $b(0)=0$ and $b$ is
absolutely continuous, $b(\rho_P(x))$ is the integral of $b'(r)$ over
$0<r<\rho_P(x)$. The same exchange of integrals as above gives
\begin{align}
 \int_D b(\rho_P(x))\dd x
 &=\int_0^\infty b'(r)|D\cap\{\rho_P>r\}|\dd r\notag\\
 &\geq\int_0^{(nv_d)^{-1/d}}b'(r)\bigl(1-(r/(nv_d)^{-1/d})^d\bigr)\dd r\notag\\
 &=\frac d{(nv_d)^{-1}}\int_0^{(nv_d)^{-1/d}}b(r)r^{d-1}\dd r\notag\\
 &=d\int_0^1b((nv_d)^{-1/d}t)t^{d-1}\dd t\notag\\
 &\geq d\,b((nv_d)^{-1/d})\int_0^1t^{d+1}\dd t
 =\frac d{d+2}b((nv_d)^{-1/d}).
 \label{eq:threeintegral}
\end{align}
The last inequality uses $b(t(nv_d)^{-1/d})\geq t^2b((nv_d)^{-1/d})$ from
Lemma~\ref{lem:profile}. The boundary terms in integration by parts
vanish because $b(0)=0$
and $1-((nv_d)^{-1/d}/(nv_d)^{-1/d})^d=0$. Proposition~\ref{prop:integralprofile} proves
\eqref{eq:exactthree}.

By Lemmas~\ref{lem:profile} and \ref{lem:volume},
\begin{align}
 \frac d{d+2}b((nv_d)^{-1/d})
 &\geq\frac d{4(d+2)}
       \min\{A,L(nv_d)^{-1/d},H(nv_d)^{-2/d}\},\label{eq:paramintermediate}\\
 (nv_d)^{-1/d}&\geq\sqrt{\frac{d+2}{2\pi\e}}\,n^{-1/d}.
 \notag
\end{align}
Treating the three entries of the minimum separately yields
\begin{align*}
 \frac d{4(d+2)}A&\geq A/12,\\
 \frac d{4(d+2)}L(nv_d)^{-1/d}
 &\geq\frac{L\sqrt d}{4\sqrt{6\pi\e}}\,n^{-1/d},
 \\
 \frac d{4(d+2)}H(nv_d)^{-2/d}
 &\geq\frac{Hd}{8\pi\e}\,n^{-2/d}.
\end{align*}
Since $1/12\geq1/(8\pi\e)$ and
$1/(4\sqrt{6\pi\e})\geq1/(8\pi\e)$, we obtain \eqref{eq:three}.

The algorithm that always returns zero has error at most $A$.
For the midpoint rule on $m^d$ cells, the Lipschitz estimate and
Cauchy--Schwarz give
\begin{align}
 |I_Df-M_mf|
 &\leq L\sum_Q\int_Q|x-c_Q|\dd x\notag\\
 &\leq L\left(\sum_Q\int_Q|x-c_Q|^2\dd x\right)^{1/2}
 =\frac{L\sqrt d}{\sqrt{12}\,m}.
 \label{eq:Lmid}
\end{align}
The Taylor argument gives
\begin{equation}
\label{eq:Hmid}
 |I_Df-M_mf|\leq\frac{Hd}{24m^2}.
\end{equation}
Choose $m=\lfloor n^{1/d}\rfloor$. Use the midpoint rule or the
zero algorithm, whichever gives the smaller stated bound. Equations \eqref{eq:Lmid}--\eqref{eq:Hmid} prove
\eqref{eq:threeupper}.
\end{proof}

\begin{corollary}\label{cor:shapeinf}
\begin{equation}
\label{eq:shapeinf}
 \begin{aligned}
 &\frac1{8\pi\e}\min\{A,L\sqrt d\,n^{-1/d},Hd\,n^{-2/d}\}\\
 &\quad\leq\inf_{\substack{D\text{ open}\\|D|=1}}e_n(I_D,\CD)
 \leq\min\{A,L\sqrt d\,n^{-1/d},Hd\,n^{-2/d}\}.
 \end{aligned}
\end{equation}
The bounds also hold for the global class $\F$: the lower bound allows
adaptive local information, and the upper bound uses point values.
\end{corollary}
\begin{proof}
The lower estimate is uniform over every $D$. The cube is an admissible
set for the infimum, and all three coefficients in
\eqref{eq:threeupper} are at most one.
\end{proof}

A lower bound uniform over domains does not give a matching upper
bound for each domain. If the domain consists of many small, widely
separated components, $n$ samples leave most components unobserved.
The error can then be arbitrarily close to the initial error $A$.

\begin{proposition}\label{prop:shapesup}
For each fixed $d,n\in\N$ and $A,L,H>0$,
\begin{equation*}
 \sup_{\substack{D\text{ open, bounded}\\|D|=1}}
 e_n(I_D,\CD)=A.
\end{equation*}
The supremum also equals $A$ for inputs in $\F$ and adaptive local information,
even when sample points for those global inputs are allowed outside $D$.
\end{proposition}
\begin{proof}
The zero algorithm gives the upper bound $A$. Choose an integer $M>n$,
and let $D$ be the disjoint union of $M$ balls of equal volume $1/M$.
Their common radius is $(Mv_d)^{-1/d}$. Choose the centers with
pairwise distance strictly greater than $2(R+(Mv_d)^{-1/d})$, where $R$ is the
saturation radius in \eqref{eq:radii}.

A ball of radius $R$ around any point can intersect at most one of the
components of $D$: if it met two components, their centers would be at
distance at most $2(R+(Mv_d)^{-1/d})$. Consequently, the union of radius-$R$ balls centered at at most $n$ sample points
meets at most $n$ components. On each of
the remaining components, $b(\rho_P)=A$. Therefore
$\int_D b(\rho_P(x))\dd x\geq A\frac{M-n}{M}$.
Proposition~\ref{prop:integralprofile} applies to every algorithm and gives
$e_n\geq A(1-n/M)$. Each such $D$ is bounded, although its diameter may
increase with $M$. Letting $M\to\infty$ proves the assertion.
\end{proof}

\section{Recovery from local information}\label{sec:operator}

For integration on $\F$, complete local information has the same
minimal error as values and gradients. We prove this for every real
linear target continuous in the uniform norm, including integration
against a signed measure.

We first prove a related estimate that also applies to recovery of the
whole function. It bounds the error above and below using the distance
to the sample points. The exact equality for scalar targets follows in
Section~\ref{sec:scalarexact}.

\subsection{Distance bounds for positive linear problems}
Integration against a finite positive measure and recovery in $L_q$ have a
common property: if $f\leq g$, then their outputs have the same order.
A linear map with this property is called \emph{positive}. Positivity
allows us to turn pointwise bounds for a function into bounds for its
output.

Here is the setting that includes both examples. Let $\B$ be the space
of bounded real Borel functions on $\R^d$, with the supremum norm, and
let $S:\B\to Y$ be a bounded positive linear map. The output space $Y$
may be $\R$, $\R^k$ with the maximum norm, or $L_q(\mu)$. More generally,
$Y$ can be a real Banach lattice: a Banach space in which the usual
operations of maximum, minimum, and absolute value are defined, and
$|u|\leq|v|$ implies $\|u\|_Y\leq\|v\|_Y$. The order is compatible
with addition and multiplication by nonnegative scalars. The properties
used below are
\[
 |Su|\leq S|u|,\qquad
 0\leq u\leq v\ \Longrightarrow\ \|Su\|_Y\leq\|Sv\|_Y.
\]

Recall that $\rho_P(x)$ is the distance from $x$ to the finite sample
set $P$. Consider
\[
 \|S(\min\{A,L\rho_P,H\rho_P^2/2\})\|_Y.
\]
For integration, this is the explicit quantity
\[
 \int_{\R^d}
       \min\{A,L\rho_P(x),H\rho_P(x)^2/2\}\dd\mu(x).
\]
The three terms are the error bounds obtained by approximating $f(x)$
by zero, by the value at a nearest sample point, or by its first-order
Taylor polynomial there.

\begin{theorem}\label{thm:operator}
For every finite nonempty sample set $P$, the minimal errors on $\F$
satisfy
\[
 \begin{aligned}
 &\tfrac12\|S(\min\{A,L\rho_P,H\rho_P^2/2\})\|_Y
 \leq\|S(b(\rho_P))\|_Y\\
 &\quad\leq e_S^{\loc}(P)\leq e_S^{\vg}(P)\\
 &\quad\leq\|S(\min\{A,L\rho_P,H\rho_P^2/2\})\|_Y.
 \end{aligned}
\]
Thus the distance estimate determines the error within a factor of two,
even if complete local information is available. With a budget of $n$
sample points,
\begin{equation}
\label{eq:adoperator}
 \begin{aligned}
 &\frac12\inf_{1\leq|P|\leq n}\|S(\min\{A,L\rho_P,H\rho_P^2/2\})\|_Y\\
 &\quad\leq e_n^{\loc}(S,\F)
 \leq e_n^{\vg,\non}(S,\F)\\
 &\quad\leq\inf_{1\leq|P|\leq n}\|S(\min\{A,L\rho_P,H\rho_P^2/2\})\|_Y.
 \end{aligned}
\end{equation}
The same bounds apply when each sample gives derivatives up to any
order $s\geq1$, or all derivatives. Sample points can be restricted
to any prescribed nonempty set, provided that the infima use the
same restriction.
\end{theorem}

\begin{proof}[Proof of the lower bounds]
Fix $P$. Theorem~\ref{thm:lift} supplies a nonnegative smooth function
$f_{\delta,\eta}$ that vanishes near $P$ and is at least
$b(\rho_P)-L(\delta+\eta)$. Applying the positive map $S$ gives
\[
 S(b(\rho_P))\leq Sf_{\delta,\eta}
                    +L(\delta+\eta)S\one.
\]
Both terms on the right are positive. Monotonicity of the norm and
the triangle inequality imply
\begin{equation}
\label{eq:normprofile}
 \|S(b(\rho_P))\|_Y\leq\|Sf_{\delta,\eta}\|_Y
             +L(\delta+\eta)\|S\one\|_Y.
\end{equation}
The functions $f_{\delta,\eta}$ and $-f_{\delta,\eta}$ give zero local
information at $P$. Lemma~\ref{lem:zero-data} therefore gives an error
of at least $\|Sf_{\delta,\eta}\|_Y$ for any algorithm whose sample
points on the zero input form $P$. Letting $\eta=\delta/2$ and
$\delta\downarrow0$ proves the lower bound
$\|S(b(\rho_P))\|_Y$. Since $b(r)\geq\tfrac12\min\{A,Lr,Hr^2/2\}$ by
\eqref{eq:bchi}, this is at least $\|S(\min\{A,L\rho_P,H\rho_P^2/2\})\|_Y/2$.

For an adaptive algorithm, take $P$ to be its sample set on the zero
input, and then take the infimum over all sets with at most $n$ points.
If the algorithm takes no samples on zero, the constant inputs $\pm A$
give error at least $A\|S\one\|_Y$. This is at least $\|S(\min\{A,L\rho_P,H\rho_P^2/2\})\|_Y$ for
every $P$, since $0\leq\min\{A,L\rho_P,H\rho_P^2/2\}\leq A$.
\end{proof}

For the upper bound we construct an approximation from values and
gradients. Enumerate $P=\{p_1,\ldots,p_k\}$. At $x$, choose a nearest
point $p\in P$, breaking ties by the least index. The three choices just described satisfy
\begin{align*}
 |f(x)|&\leq A,\\
 |f(x)-f(p)|&\leq L|x-p|,\\
 |f(x)-f(p)-\nabla f(p)\cdot(x-p)|&\leq H|x-p|^2/2.
\end{align*}
The last two estimates follow from the gradient bound and Taylor's
theorem. Define $R_Pf(x)$ by choosing, respectively, $0$, $f(p)$, or
$f(p)+\nabla f(p)\cdot(x-p)$ according to which error bound is smallest.
Fix a rule for ties. This choice depends on $x,P,A,L,H$ alone.

\begin{lemma}\label{lem:reconstruction}
The function $R_Pf$ is bounded and Borel measurable, is linear in the
values and gradients at $P$, and satisfies
\begin{equation}
\label{eq:reconstructionerror}
 |f-R_Pf|\leq\min\{A,L\rho_P,H\rho_P^2/2\}.
\end{equation}
\end{lemma}
\begin{proof}
The nearest-point regions and the regions on which each approximation
is selected are Borel sets. The three estimates above prove
\eqref{eq:reconstructionerror}. Since $|f|\leq A$ and $\min\{A,L\rho_P,H\rho_P^2/2\}\leq A$,
we also have $|R_Pf|\leq2A$. Each of the three approximations is linear
in the data, and the choice between them is independent of the data.
\end{proof}

\begin{proof}[Completion of the proof of Theorem~\ref{thm:operator}]
Apply $S$ to the approximation. Positivity
and \eqref{eq:reconstructionerror} give
\[
 \|Sf-S(R_Pf)\|_Y
 \leq\|S|f-R_Pf|\|_Y
 \leq\|S(\min\{A,L\rho_P,H\rho_P^2/2\})\|_Y.
\]
This proves the upper bound for fixed $P$, and taking the infimum over
$P$ proves the upper bound in \eqref{eq:adoperator}. Data containing all
derivatives up to an order $s\geq1$ include the value and gradient, and complete local
information determines every derivative. All the other information
models in the theorem therefore satisfy the same bounds.
\end{proof}

For simultaneous integration against positive measures
$\mu_1,\ldots,\mu_k$, with error measured in the maximum norm, the
quantity is $\max_j\int\min\{A,L\rho_P,H\rho_P^2/2\}\dd\mu_j$.
For fixed $S$ and $P$, applying $S$ to $R_Pf$ uses fixed coefficients
or fixed output functions.

\subsection{Why the factor two is necessary for fixed points}
The following one-dimensional example attains the lower factor $1/2$
in Theorem~\ref{thm:operator}. Its upper bound uses an interpolation
rule based on values and first derivatives; related rules appear in
\cite{BBS2010,BS2011,Borodachov2022,Bojanov1975}.

\begin{proposition}\label{prop:halfsharp}
Take $d=A=L=H=1$, $P=\{-1,1\}$, and $Sf=f(0)$. Then
$\|S(\min\{A,L\rho_P,H\rho_P^2/2\})\|_Y=1/2$ and $e_S^{\loc}(P)=e_S^{\vg}(P)=1/4$.
Thus the constant $1/2$ in Theorem~\ref{thm:operator} is optimal for
prescribed sample points.
\end{proposition}
\begin{proof}
Here $\rho_P(0)=1$, $\min\{1,1,1/2\}=1/2$, and $b(1)=1/4$, so
Theorem~\ref{thm:operator} gives the lower bound $1/4$.
For the upper bound, use
\[
 \mathcal Hf=\tfrac12(f(-1)+f(1))
                  +\tfrac14(f'(-1)-f'(1)).
\]
This rule is exact on affine functions. Twice integrating the second
derivative gives
\begin{equation}
\label{eq:Peano}
 f(0)-\mathcal Hf
 =\int_{-1}^1\left(\frac{|t|}{2}-\frac14\right)f''(t)\dd t.
\end{equation}
To check this identity, substitute
$f(x)=f(-1)+(x+1)f'(-1)+\int_{-1}^x(x-t)f''(t)\dd t$
and $f'(1)=f'(-1)+\int_{-1}^1f''(t)\dd t$ into the rule.
The coefficient of $f''(t)$ is $-t/2-1/4$ for $-1<t<0$ and
$t/2-1/4$ for $0<t<1$. Hence
\begin{align*}
 \int_{-1}^1\left|\frac{|t|}{2}-\frac14\right|\dd t
 &=2\left[\int_0^{1/2}(\tfrac14-\tfrac t2)\dd t
          +\int_{1/2}^1(\tfrac t2-\tfrac14)\dd t\right]\\
 &=\tfrac14.
\end{align*}
Together with $\|f''\|_\infty\leq1$, this proves the upper bound.
\end{proof}

The points must be fixed in this example: if they can be chosen freely,
one sample at zero recovers $f(0)$ exactly.

\subsection{Making a function zero near the sample points}
\label{sec:flatten}

To obtain an exact comparison, we need to relate two conditions:
$f$ and $\nabla f$ vanish at each sample point, and $f$ vanishes on a
neighborhood of each sample point. The second condition is stronger,
but it can be achieved by an arbitrarily small change in the uniform
norm without leaving $\F$. A target continuous in that norm then
converges along with the approximating functions.

\begin{theorem}\label{thm:density}
Let $P$ be finite and nonempty. If $f\in\F$ satisfies
$f(p)=0$ and $\nabla f(p)=0$ for every $p\in P$, then there are
functions $f_j\in\F$, each identically zero on a neighborhood of $P$,
such that $\|f_j-f\|_\infty\to0$.
\end{theorem}

The proof multiplies $f$ by a cutoff near each point. A cutoff changing
from zero to one across a ball of radius $a$ has second derivatives
of order $a^{-2}$. Even though $f$ is of order $a^2$ there, their
product need not be small. We instead make the cutoff change slowly
as a function of the logarithm of the distance. This reduces the
extra terms in the Hessian to order $1/T$, where $T$ can be arbitrarily
large.

Choose a smooth function $\zeta:\R\to[0,1]$ that is zero on
$(-\infty,0]$ and one on $[1,\infty)$. Fix a constant $C\geq1$ bounding
both $|\zeta'|$ and $|\zeta''|$. For $a>0$ and $T\geq1$, define
\[
 \theta_p(x)=
 \begin{cases}
 \zeta\bigl(1+T^{-1}\log(|x-p|/a)\bigr),&x\ne p,\\
 0,&x=p.
 \end{cases}
\]
This function is zero on $B(p,a\e^{-T})$ and one outside $B(p,a)$.
It is smooth, including at $p$ and at the boundaries of these balls.

\begin{lemma}\label{lem:cutoff}
For $x\ne p$,
\begin{equation}
\label{eq:cutoffhessian}
 |\nabla\theta_p(x)|\leq\frac C{T|x-p|},\qquad
 \|\nabla^2\theta_p(x)\|_\op\leq\frac{2C}{T|x-p|^2}.
\end{equation}
\end{lemma}
\begin{proof}
For the radial function $r\mapsto\zeta(1+T^{-1}\log(r/a))$, $r>0$, the derivative with respect to $r$ is
$\zeta'(1+T^{-1}\log(r/a))/(Tr)$, and the second derivative is
$\zeta''(1+T^{-1}\log(r/a))/(T^2r^2)-\zeta'(1+T^{-1}\log(r/a))/(Tr^2)$.
For a radial function, the Hessian has the second radial derivative
as its radial eigenvalue, and the first radial derivative divided by
$r$ as each tangential eigenvalue. Their absolute values are bounded
by $2C/(Tr^2)$ because $T\geq1$. In dimension one only the radial
eigenvalue occurs.
\end{proof}

\begin{proof}[Proof of Theorem~\ref{thm:density}]
Choose $a\leq L/H$ small enough that the balls $B(p,a)$ are pairwise
disjoint. Define $g=\theta_pf$ in $B(p,a)$ for each $p\in P$, and
$g=f$ outside these balls. Then $g$ is smooth and $|g|\leq A$.
Only one cutoff acts at each point.

For $|x-p|\leq a$, Taylor's theorem and the assumptions at $p$ give
$|f(x)|\leq H|x-p|^2/2$ and $|\nabla f(x)|\leq H|x-p|$. Thus
\begin{align}
 |\nabla g(x)|
 &\leq\theta_p(x)|\nabla f(x)|
      +|f(x)|\,|\nabla\theta_p(x)|\notag\\
 &\leq L+\frac{H|x-p|^2}{2}\frac C{T|x-p|}
 \leq L\left(1+\frac C{2T}\right).
 \label{eq:cutoffgradientbound}
\end{align}
The product rule for the Hessian gives
\[
 \nabla^2g=\theta_p\nabla^2f
    +\nabla\theta_p\otimes\nabla f
    +\nabla f\otimes\nabla\theta_p
    +f\nabla^2\theta_p.
\]
Here $u\otimes v$ is the matrix $uv^{\mathsf T}$, whose operator norm
is $|u|\,|v|$. Therefore
\begin{align}
 \|\nabla^2g(x)\|_\op
 &\leq H+2\frac C{T|x-p|}H|x-p|
      +\frac{H|x-p|^2}{2}\frac{2C}{T|x-p|^2}\notag\\
 &\leq H\left(1+\frac{3C}{T}\right).
 \label{eq:cutoffHbound}
\end{align}
Outside the balls, $g=f$, so these estimates hold everywhere.

Divide $g$ by $1+3C/T$. The
function, gradient, and Hessian bounds now give $g/(1+3C/T)\in\F$. It vanishes
on every $B(p,a\e^{-T})$, and
\begin{align}
 \left\|\frac{g}{1+3C/T}-f\right\|_\infty
 &\leq\left\|\frac{g}{1+3C/T}-g\right\|_\infty+\|g-f\|_\infty\notag\\
 &\leq\frac{3AC}{T}+\frac{Ha^2}{2}.
 \label{eq:densityerror}
\end{align}
Let $T\to\infty$ and $a\downarrow0$, keeping the balls disjoint.
Choosing sequences of these parameters gives the required $f_j$.
\end{proof}

Only uniform convergence is asserted. For example, if
$\nabla^2f(p)\ne0$, convergence of the Hessians is impossible because
each $f_j$ is zero near $p$.

\subsection{Exact errors for scalar targets}\label{sec:scalarexact}
Let $S$ be a real linear functional on the smooth functions with
bounded function, gradient, and Hessian. Assume that $S$ is continuous
in the uniform norm: there is a finite constant $C_S$ such that, for
every such $f$,
\begin{equation}
\label{eq:Ssupcontinuous}
 |Sf|\leq C_S\|f\|_\infty.
\end{equation}
Integration against a finite signed measure satisfies this condition.

For a finite nonempty set $P$, define
\begin{equation}
\label{eq:radiusP}
 r_P(S)=\sup\bigl\{|Sf|:f\in\F,
           \ f(p)=0,\ \nabla f(p)=0\text{ for every }p\in P\bigr\}.
\end{equation}
This is the largest absolute target value that can remain when all observed
values and gradients are zero. The next lemma identifies it with the
minimal error from values and gradients. It follows from the
Hahn--Banach theorem; see \cite{Bakhvalov1971,NW2008}.

\begin{lemma}\label{lem:HB}
For each finite nonempty $P$, one has $e_S^{\vg}(P)=r_P(S)$.
A linear combination of the observed values and gradient components
attains this error.
\end{lemma}
\begin{proof}
For this proof, give the space $X$ of smooth functions with bounded
function, gradient, and Hessian the norm
\begin{equation}
\label{eq:Xnorm}
 \|f\|_{A,L,H}
 =\max\left\{\frac{\|f\|_\infty}{A},
              \frac{\|\nabla f\|_\infty}{L},
              \frac{\|\nabla^2f\|_{\op,\infty}}{H}\right\}.
\end{equation}
Its unit ball is $\F$. Let $Nf$ be the vector of values and gradients
at the points in $P$. On the subspace $\ker N$, the norm of $S$ is
$r_P(S)$. Hahn--Banach extends this restriction to a linear functional
$T$ on $X$ with the same norm.

Since $S-T$ vanishes on $\ker N$, any two functions with the same
data have the same value of $S-T$. Thus $S-T$ is a linear function
of $Nf$: write $S-T=\ell\circ N$, extending $\ell$ to the full
finite-dimensional data space if necessary. The rule
$\ell(Nf)$ has error
\[
 \sup_{f\in\F}|Sf-\ell(Nf)|
 =\sup_{f\in\F}|Tf|=r_P(S).
\]
For the converse, $f$ and $-f$ give identical data whenever $Nf=0$.
Lemma~\ref{lem:zero-data} gives error at least $|Sf|$; taking the
supremum proves the lower bound.
\end{proof}

\begin{theorem}\label{thm:exactlocal}
Let $S$ satisfy \eqref{eq:Ssupcontinuous} and let $n\in\N$. On $\F$,
adaptive local information has the same minimal error as values and
gradients at points chosen in advance:
\begin{align}
 e_n^{\loc}(S,\F)
 &=e_n^{\vg}(S,\F)
  =e_n^{\vg,\non}(S,\F)\notag\\
 &=\inf_{1\leq|P|\leq n}r_P(S).
 \label{eq:exactlocal}
\end{align}
For each fixed finite nonempty $P$,
\begin{equation}
\label{eq:fixedexact}
 e_S^{\loc}(P)=e_S^{\vg}(P)=r_P(S).
\end{equation}
The same errors result from observing derivatives up to any order
$s\geq1$, or all derivatives, with adaptive or nonadaptive points.
All comparisons remain valid when the sample points are restricted
to the same prescribed nonempty set.
\end{theorem}
\begin{proof}
First, Theorem~\ref{thm:density} and uniform-norm continuity show that
\[
 \sup_{\substack{f\in\F\\ f=0\text{ near }P}}|Sf|=r_P(S).
\]
The left side is at most the right side, since vanishing near $P$
implies that the values and gradients vanish there. Conversely, if
$f\in\F$ has zero values and gradients on $P$, the theorem gives
$f_j\in\F$ vanishing near $P$ with
$|Sf_j-Sf|\leq C_S\|f_j-f\|_\infty\to0$. This proves equality.

Run an adaptive algorithm with local information on zero and denote
its sample set by $P_0$. If it takes no samples, choose any one
admissible point for $P_0$. Every function vanishing near $P_0$ gives
the same observations as zero. The preceding equality and
Lemma~\ref{lem:zero-data} imply
\begin{equation}
\label{eq:adexactlower}
 e(\mathcal A,S,\F)
 \geq r_{P_0}(S)
 \geq\inf_{1\leq|P|\leq n}r_P(S).
\end{equation}
For each admissible $P$, Lemma~\ref{lem:HB} gives a rule with fixed
points, using values and gradients, whose error is $r_P(S)$. Hence
\begin{equation}
\label{eq:adexactupper}
 e_n^{\vg,\non}(S,\F)\leq\inf_{1\leq|P|\leq n}r_P(S).
\end{equation}
Every information model in the theorem is at least as powerful as
this fixed-point rule and no more powerful than adaptive local
information. The lower and upper bounds therefore coincide for all
of them. The same argument with $P$ fixed proves
\eqref{eq:fixedexact}.
\end{proof}

\begin{corollary}\label{cor:signed}
For every finite signed Borel measure $\nu$ on $\R^d$, integration
$S=I_\nu$ satisfies \eqref{eq:exactlocal}, since
$|I_\nu f|\leq |\nu|(\R^d)\|f\|_\infty$.
For positive $\nu$, the common error also satisfies
\eqref{eq:adoperator}.
\end{corollary}

The continuity assumption distinguishes integration from derivative
recovery. In one dimension, for example, fix the sample set $P=\{p\}$. The
target $Sf=f''(p)$ is determined by second-derivative data at $p$,
but is not determined by $f(p)$ and $f'(p)$. It does not satisfy \eqref{eq:Ssupcontinuous}.

\section{Function recovery and other integration measures}\label{sec:extensions}
\subsection{Recovery in \texorpdfstring{$L_q$}{Lq}}\label{sec:Lq}

Suppose that the output is an approximation to the function itself.
For a finite positive Borel measure $\mu$, we measure its error by
\[
 \|f-\widehat f\|_{L_q(\mu)}
 =\left(\int |f-\widehat f|^q\dd\mu\right)^{1/q}
 \qquad(1\leq q<\infty).
\]
For $q=\infty$, the error is the essential supremum of
$|f-\widehat f|$: a bound may fail on a set of measure zero.
Write $e_{n,q}^{\mathcal I}$ for the minimal worst-case error using
at most $n$ sample points and information of type $\mathcal I$.
The identity map from bounded Borel functions to $L_q(\mu)$ is a bounded
positive linear operator. Applying Theorem~\ref{thm:operator} gives
\begin{equation}
\label{eq:Lqquant}
 \begin{aligned}
 &\frac12\inf_{1\leq|P|\leq n}\|\min\{A,L\rho_P,H\rho_P^2/2\}\|_{L_q(\mu)}\\
 &\quad\leq e_{n,q}^{\loc}(\F,\mu)
 \leq e_{n,q}^{\vg,\non}(\F,\mu)\\
 &\quad\leq\inf_{1\leq|P|\leq n}\|\min\{A,L\rho_P,H\rho_P^2/2\}\|_{L_q(\mu)}.
 \end{aligned}
\end{equation}
For Lebesgue measure on a set of volume one, the distance estimate also
gives explicit constants. Here the lower bound uses the function $b$
itself, before the comparison $b(r)\geq\tfrac12\min\{A,Lr,Hr^2/2\}$.

\begin{theorem}\label{thm:Lq}
Let $D\subset\R^d$ be measurable with $|D|=1$.
For recovery on $D$ of functions in $\Fn$, the minimal errors satisfy,
for $1\leq q<\infty$,
\begin{equation}
\label{eq:Lqexplicit}
 e_{n,q}^{\loc}
 \geq\frac{v_d^{-2/d}}{4d}
       \left(\frac d{d+2q}\right)^{1/q}n^{-2/d}.
\end{equation}
For $q=\infty$,
\begin{equation}
\label{eq:Linfexplicit}
 e_{n,\infty}^{\loc}
 \geq\frac{v_d^{-2/d}}{4d}n^{-2/d}.
\end{equation}
In particular, uniformly for $1\leq q\leq\infty$,
\begin{equation}
\label{eq:Lquniform}
 e_{n,q}^{\loc}\geq\frac1{8\pi\e}n^{-2/d}.
\end{equation}
For open $D$, these bounds also hold on $\widetilde C^2(D)$ when
algorithms use local information at points of $D$.
\end{theorem}
\begin{proof}
Let $P$ be the points sampled by the algorithm on the zero input.
The smooth functions from Theorem~\ref{thm:lift} vanish near $P$.
They and their negatives therefore give the same data as zero.
Lemma~\ref{lem:zero-data}, followed by the limit in
\eqref{eq:liftprofile}, gives
$e(\mathcal A)\geq\|b(\rho_P)\|_{L_q(D)}$.
Recall that $b(r)=r^2/(4d)$ for
$0\leq r\leq (nv_d)^{-1/d}$. For $q<\infty$, integrate the derivative of
$(r^2/(4d))^q$ and apply Tonelli's theorem and
\eqref{eq:tubevolume}:
\begin{align*}
 \int_D b(\rho_P(x))^q\dd x
 &\geq\frac{2q}{(4d)^q}
      \int_0^{(nv_d)^{-1/d}}r^{2q-1}(1-nv_dr^d)\dd r\\
 &=\frac1{(4d)^q}
     \left((nv_d)^{-2q/d}-\frac{2q}{d+2q}(nv_d)^{-2q/d}\right)\\
 &=\frac{(nv_d)^{-2q/d}}{(4d)^q}\frac d{d+2q}.
\end{align*}
Taking $q$th roots proves \eqref{eq:Lqexplicit}.
For $q=\infty$, we need a set of positive volume on which the function
is large. For every $0<r<(nv_d)^{-1/d}$, \eqref{eq:tubevolume} gives
$|D\cap\{\rho_P>r\}|>0$. Since $b$ is nondecreasing,
$\|b(\rho_P)\|_{L_\infty(D)}\geq b(r)$.
Let $r\uparrow (nv_d)^{-1/d}$ to obtain \eqref{eq:Linfexplicit}.
Finally, on a probability space
$\|u\|_{L_q(D)}\geq\|u\|_{L_1(D)} \qquad(1\leq q\leq\infty)$.
Use \eqref{eq:mainintegral} and \eqref{eq:upperball} to conclude
\eqref{eq:Lquniform}.
\end{proof}

On the unit cube, divide each coordinate interval into $m$ equal parts
and use the first-order Taylor polynomial at each of the $m^d$ cell
centers. This gives a piecewise affine approximation. On a cell $Q$ with
center $c_Q$, the Hessian bound gives
\begin{equation*}
 |f(x)-f(c_Q)-\nabla f(c_Q)\cdot(x-c_Q)|
 \leq\frac1{2d}|x-c_Q|^2
 \leq\frac1{8m^2}.
\end{equation*}
Choose $m=\lfloor n^{1/d}\rfloor\geq n^{1/d}/2$. The pointwise bound
also bounds every $L_q$ norm on the unit cube, so
\begin{equation*}
 e_{n,q}^{\vg,\non}\bigl(\Fn,(0,1)^d\bigr)
 \leq\frac12n^{-2/d},\qquad1\leq q\leq\infty.
\end{equation*}

\subsection{Integration when the measure of small balls is bounded}\label{sec:measures}

The volume argument also applies to weighted integration
$I_\mu f=\int f\dd\mu$. Let $\mu$ be a finite positive Borel measure
on $\R^d$, with total mass $m=\mu(\R^d)>0$.
We now bound how much of this mass can lie near the sample points.
For $r\geq0$, let
$\beta_\mu(r)=\sup_{z\in\R^d}\mu(\overline B(z,r))$.
Thus $\beta_\mu(r)$ is the largest mass of a closed ball of radius $r$.
The set $\{\rho_P\leq r\}$ is the union of the balls centered at $P$.
For $|P|\leq n$, adding their masses gives
$\mu\{\rho_P\leq r\}\leq n\beta_\mu(r)$ and hence
$\mu\{\rho_P>r\}\geq(m-n\beta_\mu(r))_+$.

\begin{theorem}\label{thm:smallball}
For $\F$ and deterministic adaptive local information,
\begin{equation}
\label{eq:smallballbound}
  e_n^{\loc}(I_\mu,\F)\geq
 \int_0^\infty b'(r)(m-n\beta_\mu(r))_+\dd r.
\end{equation}
Suppose that, for some $B,s>0$,
\begin{equation}
\label{eq:powerball}
 \beta_\mu(r)\leq Br^s\qquad(r>0).
\end{equation}
At radius $(m/(nB))^{1/s}$, one has $nB(m/(nB))=m$. Then
\begin{equation*}
 e_n^{\loc}(I_\mu,\F)
 \geq ms\int_0^1b((m/(nB))^{1/s}t)t^{s-1}\dd t
 \geq\frac{ms}{s+2}b((m/(nB))^{1/s}).
\end{equation*}
\end{theorem}
\begin{proof}
Let $P$ be the points sampled on the zero input. Integrating the
uniform lower estimate in Theorem~\ref{thm:lift} and letting
$\delta,\eta\downarrow0$ gives
$e(\mathcal A)\geq\int b(\rho_P)\dd\mu$.
Tonelli's theorem gives
\begin{align*}
 \int b(\rho_P(x))\dd\mu(x)
 &=\int_0^\infty b'(r)\mu\{\rho_P>r\}\dd r\\
 &\geq\int_0^\infty b'(r)(m-n\beta_\mu(r))_+\dd r.
\end{align*}
Under \eqref{eq:powerball}, $m-nBr^s$ is nonnegative for $0\leq r\leq(m/(nB))^{1/s}$.
Integrate over this interval, use integration by parts, and then apply
$b((m/(nB))^{1/s}t)\geq t^2b((m/(nB))^{1/s})$ from \eqref{eq:bscale}:
\begin{align}
 \int_0^{(m/(nB))^{1/s}} b'(r)(m-nBr^s)\dd r
 &=nBs\int_0^{(m/(nB))^{1/s}}b(r)r^{s-1}\dd r\notag\\
 &=ms\int_0^1b((m/(nB))^{1/s}t)t^{s-1}\dd t\notag\\
 &\geq ms\,b((m/(nB))^{1/s})\int_0^1t^{s+1}\dd t
 =\frac{ms}{s+2}b((m/(nB))^{1/s}).
 \label{eq:smallballpowerproof}
\end{align}
The zero boundary terms use $b(0)=0$ and $nB(m/(nB))=m$.
\end{proof}

The full evaluation of the integral in
\eqref{eq:smallballpowerproof} is given in Appendix~\ref{app:measure-integral}.
If
\begin{equation}
\label{eq:quadraticrange}
 (m/(nB))^{1/s}\leq\min\{2L/H,2\sqrt{A/H}\},
\end{equation}
then $b(r)=Hr^2/4$ for $0\leq r\leq(m/(nB))^{1/s}$, so
\begin{equation}
\label{eq:powerquadratic}
 e_n^{\loc}(I_\mu,\F)
 \geq\frac{Hm}{4}\frac{s}{s+2}
                 \left(\frac m{nB}\right)^{2/s}.
\end{equation}
The exponent $2/s$ comes from the growth of the measure of balls.
For volume measure restricted to a finite-volume set, one may take
$s=d$; the hypothesis also allows noninteger $s$.

\subsubsection{Bounded densities}
If $\mu$ has a Lebesgue density $w$ with $0\leq w\leq M$, then
$\beta_\mu(r)\leq Mv_dr^d$.
Theorem~\ref{thm:smallball} applies with $B=Mv_d$ and $s=d$.
If $(m/(nMv_d))^{1/d}$ satisfies \eqref{eq:quadraticrange},
\begin{equation}
\label{eq:densitylower}
 e_n^{\loc}(I_\mu,\F)
 \geq\frac{Hm}{4}\frac d{d+2}
                  \left(\frac m{nMv_d}\right)^{2/d}.
\end{equation}
Without the restriction \eqref{eq:quadraticrange},
\begin{equation*}
 e_n^{\loc}(I_\mu,\F)
 \geq\frac{md}{d+2}\,
 b\left(\left(\frac m{nMv_d}\right)^{1/d}\right).
\end{equation*}

\subsubsection{Standard Gaussian integration}
Let $\gamma_d$ be the standard Gaussian probability measure:
$\dd\gamma_d(x)=(2\pi)^{-d/2}\e^{-|x|^2/2}\dd x$.

\begin{corollary}\label{cor:gaussian}
For every $d,n\in\N$,
\begin{equation*}
 e_n^{\loc}(I_{\gamma_d},\Fn)
 \geq\frac{2\pi v_d^{-2/d}}{4(d+2)}n^{-2/d}
 \geq\frac1{4\e}n^{-2/d}.
\end{equation*}
\end{corollary}
\begin{proof}
Here $m=1$ and $M=(2\pi)^{-d/2}$. The relevant radius is
$\sqrt{2\pi}\,v_d^{-1/d}n^{-1/d} \leq\sqrt{\pi d/2}\,n^{-1/d}<2\sqrt d$.
The inequality uses \eqref{eq:lowerball}. For the normalized profile,
both thresholds in \eqref{eq:quadraticrange} equal $2\sqrt d$.
Therefore \eqref{eq:densitylower} applies with $H=1/d$ and gives the
first estimate. Equation \eqref{eq:upperball} yields the second.
\end{proof}

For measures with atoms, use \eqref{eq:smallballbound} directly.
Indeed, if $\mu(\{x\})>0$, then $\beta_\mu(r)\geq\mu(\{x\})$ for every
$r>0$, whereas $Br^s\to0$ as $r\downarrow0$ when $s>0$.

\subsection{A change of metric and arbitrary volume}\label{sec:anisotropic}

Different smoothness bounds in different directions can be handled
by a linear change of coordinates. Let $G$ be a real symmetric positive
definite $d\times d$ matrix, and let $G^{1/2}$ be its positive definite
square root. Distances after the change of coordinates are
$|x|_G=|G^{1/2}x|$. For example, if $G$ is diagonal, this change stretches
the coordinate axes by the square roots of its diagonal entries.
Write $\rho_{G,P}(x)=\min_{p\in P}|x-p|_G$ and consider the class
\begin{equation*}
 \mathcal F_G(A,L,H)=\left\{f\in C^\infty(\R^d):
 \begin{aligned}
  &\|f\|_\infty\leq A,\quad
   \|G^{-1/2}\nabla f\|_\infty\leq L,\\
  &\|G^{-1/2}\nabla^2 f\,G^{-1/2}\|_{\op,\infty}\leq H
 \end{aligned}\right\}.
\end{equation*}

\begin{theorem}\label{thm:G}
Let $D\subset\R^d$ be measurable with volume $0<V<\infty$.
For adaptive local information on $\mathcal F_G(A,L,H)$,
\begin{align}
 e_n^{\loc}(I_D,\mathcal F_G)
 &\geq\frac{Vd}{d+2}b((V\sqrt{\det G}/(nv_d))^{1/d})\label{eq:Gexact}\\
 &\geq\frac V{8\pi\e}
 \min\left\{A,
 L\sqrt d\left(\frac{V\sqrt{\det G}}n\right)^{1/d},
 Hd\left(\frac{V\sqrt{\det G}}n\right)^{2/d}\right\}.
 \label{eq:Guniform}
\end{align}
On $D=G^{-1/2}(0,(V\sqrt{\det G})^{1/d})^d$, a midpoint
rule attains the same dependence on $A,L,H,V,G,d,n$, up to absolute
constant factors.
\end{theorem}
\begin{proof}
Use the change of variables
$z=G^{1/2}x$ and $g(z)=f(G^{-1/2}z)$.
The chain rule and Jacobian give
\begin{align*}
 \nabla g(z)&=G^{-1/2}\nabla f(G^{-1/2}z),\\
 \nabla^2g(z)&=G^{-1/2}\nabla^2 f(G^{-1/2}z)G^{-1/2},
 \\
 |G^{1/2}D|&=V\sqrt{\det G},\qquad
 I_Df=(\det G)^{-1/2}I_{G^{1/2}D}g.
\end{align*}
Thus $f\in\mathcal F_G$ exactly when $g\in\F$.
Each sample point $x$ corresponds to one point $z=G^{1/2}x$.
Two functions agree near $x$ exactly when their transformed functions
agree near $z$, so the change preserves local information and sample
cost. A ball of radius $r$ in the $G$-metric has Euclidean volume
$v_dr^d/\sqrt{\det G}$. Hence
\begin{equation}
\label{eq:Gtube}
 |D\cap\{\rho_{G,P}\leq r\}|
 \leq\frac{nv_dr^d}{\sqrt{\det G}}.
\end{equation}
Apply Theorem~\ref{thm:lift} in the $z$ coordinates and return to the
$x$ coordinates. The resulting functions vanish near the sample points and are bounded
below by $b(\rho_{G,P})$ minus a constant tending to zero.
Using these functions in the lower-bound argument and integrating
\eqref{eq:Gtube} as in \eqref{eq:threeintegral} gives
\begin{align*}
 \int_Db(\rho_{G,P}(x))\dd x
 &\geq Vd\int_0^1b((V\sqrt{\det G}/(nv_d))^{1/d}t)t^{d-1}\dd t\\
 &\geq\frac{Vd}{d+2}b((V\sqrt{\det G}/(nv_d))^{1/d}),
\end{align*}
which proves \eqref{eq:Gexact}. Apply \eqref{eq:bpsi} and \eqref{eq:upperball}, using the same
coefficient comparisons as in the proof of Theorem~\ref{thm:three}, with
$n^{-1/d}$ replaced by $(V\sqrt{\det G}/n)^{1/d}$.
This proves \eqref{eq:Guniform}.

For the stated domain, work with the midpoint grid in the transformed
cube $(0,(V\sqrt{\det G})^{1/d})^d$. Its mean squared distance to a cell center is
$d(V\sqrt{\det G})^{2/d}/(12m^2)$. Transforming the amplitude, gradient, and Hessian
estimates from Theorem~\ref{thm:three}, and choosing
$m=\lfloor n^{1/d}\rfloor$, gives
\begin{equation*}
 e_n^{\val,\non}(I_D,\mathcal F_G)
 \leq V\min\left\{A,
 \frac{L\sqrt d\,(V\sqrt{\det G})^{1/d}}{\sqrt3}n^{-1/d},
 \frac{Hd(V\sqrt{\det G})^{2/d}}{6}n^{-2/d}\right\}.
\end{equation*}
Together with \eqref{eq:Guniform}, this proves the claimed order.
\end{proof}

When $G=I$, Theorem~\ref{thm:G} gives the arbitrary-volume version
\begin{equation*}
 e_n^{\loc}(I_D,\F)\geq\frac V{8\pi\e}
 \min\{A,L\sqrt d\,(V/n)^{1/d},Hd\,(V/n)^{2/d}\}.
\end{equation*}
Setting $V=1$ also recovers \eqref{eq:three} on the global class.

\section{Scope and open questions}\label{sec:scope}

The class \eqref{eq:globalclass} bounds the function and its first two
derivatives. As the smoothing radius decreases, derivatives of order
three and above may become unbounded. The constructed functions therefore
need not satisfy prescribed bounds on those derivatives.

For randomized algorithms, a lower bound on the expected squared error
must use one input independently of the random seed. The test function
in our proof may depend on that seed, so a randomized lower bound
requires another argument.

The lower bound and the cube rule determine the rate $n^{-2/d}$ up to
an absolute factor. Two constant problems remain: the best universal
constant in \eqref{eq:conjecture}, and the exact infimum over domains
in \eqref{eq:shapeinf}. The factor $1/2$ in
Proposition~\ref{prop:halfsharp} is sharp for the comparison at prescribed
sample points; optimizing the integration error over point sets and
domains is a separate problem. Further questions about asymptotic
constants and other smoothness norms appear in
\cite{Novak2020,Novak2024}.

\appendix
\section{The smooth-function construction}\label{app:construction}

\subsection{Properties of the distance function}\label{app:profile}
\begin{proof}[Proof of Lemma~\ref{lem:profile}]
Before truncation,
\begin{equation}
\label{eq:qderiv}
 q'(r)=\min\{Hr/2,L\},\qquad
 q''(r)=\begin{cases}H/2,&0<r<2L/H,\\0,&r>2L/H.
 \end{cases}
\end{equation}
At $2L/H$, both formulas for $q$ give the value $L^2/H$ and derivative
$L$. For $r\geq 2L/H$,
\[
 \frac H4r^2-q(r)=\frac H4(r-2L/H)^2\geq0,
\]
which proves \eqref{eq:bquad}. If $r\leq 2L/H$, then $Hr^2/2\leq Lr$ and
\[
 b(r)=\min\{A,Hr^2/4\}\geq\tfrac12\min\{A,Hr^2/2\}.
\]
If $r\geq 2L/H$, then $Hr^2/2\geq Lr$ and
$q(r)=Lr-L^2/H\geq Lr/2$. Together with $q\leq Lr$ and
\eqref{eq:bquad}, these inequalities prove \eqref{eq:bchi}.
Since $\min\{A,Lr,Hr^2/2\}\geq\tfrac12\min\{A,Lr,Hr^2\}$, they also give
\eqref{eq:bpsi}.

For positive $r$, the ratio $q(r)/r^2$ equals $H/4$ on the first
interval. On the second interval,
\[
 \frac{\dd}{\dd r}\frac{q(r)}{r^2}
 =\frac{-Lr+2L^2/H}{r^3}\leq0.
\]
Thus $b(r)/r^2=\min\{A/r^2,q(r)/r^2\}$ is nonincreasing, which proves
\eqref{eq:bscale} for $r,t>0$. The cases $r=0$ and $t=0$ follow from
$b(0)=0$. The other assertions follow from \eqref{eq:qderiv} and
$q(R)=A$.
\end{proof}

\subsection{Proof of the smooth construction}
\subsubsection{Making the function zero near \texorpdfstring{$K$}{K}}
The construction has three steps. First, replace the distance by
$(\rho_K-\delta)_+$ to make the function zero wherever the distance
to $K$ is at most $\delta$. Second, modify the result to control its
second differences. Third, average over a ball of radius $\eta$ to
obtain a smooth function. Since $\eta<\delta$, this last step still
leaves a zero neighborhood of $K$.

We begin with an upper bound on second derivatives. For a smooth function,
concavity after subtracting $H|x|^2/4$ says that every second directional
derivative is at most $H/2$. The same concavity condition is meaningful
for the nonsmooth function used here.

\begin{lemma}\label{lem:buffer}
The function
$h_\delta(x)=b\bigl((\rho_K(x)-\delta)_+\bigr)$
is $L$-Lipschitz, takes values in $[0,A]$, and satisfies
\begin{equation}
\label{eq:hsemi}
 x\longmapsto h_\delta(x)-\frac H4|x|^2
 \quad\text{is concave}.
\end{equation}
For every $x$ with $\rho_K(x)\leq\delta$ and every $y\in\R^d$,
\begin{equation}
\label{eq:hzeroestimate}
 h_\delta(y)\leq\frac H4|x-y|^2.
\end{equation}
\end{lemma}
\begin{proof}
For a single point $p$, consider
$u_p(x)=q\bigl((|x-p|-\delta)_+\bigr)$. At distance
$r=|x-p|>\delta$, its radial derivatives away from the two joining radii
satisfy
\[
 u_p'(r)=\min\{H(r-\delta)/2,L\},\qquad
 0\leq u_p''(r)\leq H/2.
\]
For a radial function, the Hessian has eigenvalue $u_p''(r)$ in the
radial direction and eigenvalue $u_p'(r)/r$ in every tangential direction.
The latter satisfies
\[
 \frac{u_p'(r)}r\leq\frac{H(r-\delta)}{2r}\leq H/2.
\]
Thus all second directional derivatives are at most $H/2$ wherever these
formulas apply.
For $r\leq\delta$, the function is zero. Its first derivatives are
continuous at $r=\delta$ and $r=\delta+2L/H$. On each straight line,
the first derivative is therefore continuous and piecewise smooth, with
second derivative at most $H/2$. Subtracting $H|x|^2/4$ makes the
restriction to each line concave. Hence
$u_p(x)-H|x|^2/4$ is concave on $\R^d$.

By monotonicity and continuity of $q$,
\begin{align}
 h_\delta(x)-H|x|^2/4
 =\inf\bigl(&\{A-H|x|^2/4\}\notag\\
 &\cup\{u_p(x)-H|x|^2/4:p\in K\}\bigr).
 \label{eq:infconcave}
\end{align}
Every function inside the infimum is concave. Their pointwise infimum
is concave as well: each concavity inequality is bounded below by the
same convex combination of the two endpoint infima. This proves
\eqref{eq:hsemi}. The size bound follows from $b\leq A$.
The Lipschitz bound follows because $b$ is $L$-Lipschitz and both the
distance and the positive-part function are $1$-Lipschitz.
Finally, if $\rho_K(x)\leq\delta$, the triangle inequality gives
$(\rho_K(y)-\delta)_+\leq|y-x|$. Apply \eqref{eq:bquad} to obtain
\eqref{eq:hzeroestimate}.
\end{proof}

\subsubsection{Quadratic regularization}
Define
\begin{equation}
\label{eq:envelope}
 F_\delta(x)=\sup_{y\in\R^d}
 \left\{h_\delta(y)-\frac H2|x-y|^2\right\}.
\end{equation}
For each $x$, we compare the values $h_\delta(y)$ after subtracting a
quadratic penalty for the distance from $x$ to $y$. The choice $y=x$
shows that $F_\delta$ lies above $h_\delta$. The one-sided bound in Lemma~\ref{lem:buffer} makes the new function
satisfy two-sided bounds on second differences. These bounds are expressed below by convexity
and concavity, without assuming differentiability.

\begin{lemma}\label{lem:envelope}
The supremum in \eqref{eq:envelope} is attained, and
\begin{align}
 h_\delta&\leq F_\delta\leq A,\notag\\
 \Lip(F_\delta)&\leq L,\label{eq:Flip}\\
 F_\delta(x)&=0\quad\text{if }\rho_K(x)\leq\delta,\label{eq:Fzero}\\
 F_\delta(x)+\tfrac H2|x|^2&\quad\text{is convex},\notag\\
 F_\delta(x)-\tfrac H2|x|^2&\quad\text{is concave}.\label{eq:Fconcave}
\end{align}
Equivalently, the symmetric second difference satisfies
\begin{equation}
\label{eq:seconddifference}
 -H|z|^2\leq F_\delta(x+z)+F_\delta(x-z)-2F_\delta(x)
 \leq H|z|^2.
\end{equation}
\end{lemma}
\begin{proof}
Taking $y=x$ gives $F_\delta(x)\geq h_\delta(x)\geq0$.
If $|y-x|>\sqrt{2A/H}$, the expression in braces is negative.
It is therefore enough to maximize a continuous function on a compact
ball. The bound $F_\delta\leq A$ follows from $h_\delta\leq A$.

For $x,x'\in\R^d$, use the point $y+x'-x$ in the supremum defining
$F_\delta(x')$:
\begin{align*}
 F_\delta(x')
 &\geq h_\delta(y+x'-x)-\tfrac H2|x-y|^2\\
 &\geq h_\delta(y)-\tfrac H2|x-y|^2-L|x'-x|.
\end{align*}
Take the supremum over $y$ and interchange $x,x'$ to get
\eqref{eq:Flip}. If $\rho_K(x)\leq\delta$, Lemma~\ref{lem:buffer} gives
\[
 h_\delta(y)-\tfrac H2|x-y|^2
 \leq-\tfrac H4|x-y|^2\leq0.
\]
The choice $y=x$ gives equality, proving \eqref{eq:Fzero}.

To prove convexity, expand the square:
\[
 F_\delta(x)+\tfrac H2|x|^2
 =\sup_y\{h_\delta(y)-\tfrac H2|y|^2+H\langle x,y\rangle\}.
\]
The right side is a supremum of affine functions of $x$.
For the concavity assertion, both variables matter. Rewrite the expression
inside the supremum as
\begin{align*}
 &h_\delta(y)-\tfrac H2|x-y|^2-\tfrac H2|x|^2\\
 &\hspace{1cm}=
 \bigl(h_\delta(y)-\tfrac H4|y|^2\bigr)-H|x-y/2|^2.
\end{align*}
By \eqref{eq:hsemi}, this expression is jointly concave in $(x,y)$.
Choose maximizing points $y_0,y_1$ for $x_0,x_1$ and let $0\leq t\leq1$.
Evaluating at $(1-t)y_0+ty_1$ and using joint concavity gives
\begin{align*}
 &F_\delta((1-t)x_0+tx_1)-\tfrac H2|(1-t)x_0+tx_1|^2\\
 &\quad\geq(1-t)\bigl(F_\delta(x_0)-\tfrac H2|x_0|^2\bigr)\\
 &\qquad+t\bigl(F_\delta(x_1)-\tfrac H2|x_1|^2\bigr).
\end{align*}
This proves \eqref{eq:Fconcave}. Apply the convexity and concavity
inequalities at the midpoint of $x-z$ and $x+z$ to obtain
\eqref{eq:seconddifference}.
\end{proof}

\subsubsection{Smoothing by convolution}
Choose a nonnegative $\varphi\in C_c^\infty(\R^d)$ with integral one
and support in $\{z:|z|<1\}$. Put
\[
 f_{\delta,\eta}(x)
 =\int_{\R^d}F_\delta(x-z)\eta^{-d}\varphi(z/\eta)\dd z.
\]
This convolution is a weighted average of values within distance
$\eta$ of $x$. Because the weights are nonnegative and integrate to
one, the size, Lipschitz, and second-difference bounds survive averaging. We verify the derivative bounds using this observation.

\begin{proof}[Proof of Theorem~\ref{thm:lift}]
The convolution is smooth on $\R^d$, and
\[
 0\leq f_{\delta,\eta}\leq A,\qquad
 |f_{\delta,\eta}(x)-f_{\delta,\eta}(x')|\leq L|x-x'|.
\]
Thus $\|\nabla f_{\delta,\eta}\|_\infty\leq L$.
Convolving \eqref{eq:seconddifference} gives
\[
 \left|f_{\delta,\eta}(x+z)+f_{\delta,\eta}(x-z)
 -2f_{\delta,\eta}(x)\right|\leq H|z|^2.
\]
Set $z=t\theta$ for a unit vector $\theta$, divide by $t^2$, and let
$t\to0$. Then
$|\theta^T\nabla^2f_{\delta,\eta}(x)\theta|\leq H$.
For a symmetric matrix, the operator norm equals the largest absolute
value of its quadratic form on unit vectors. Hence
$\|\nabla^2f_{\delta,\eta}\|_{\op,\infty}\leq H$.

If $\rho_K(x)\leq\delta-\eta$ and $|z|\leq\eta$, then
$\rho_K(x-z)\leq\delta$. Equation~\eqref{eq:Fzero} therefore proves
\eqref{eq:liftflat}. For $|z|\leq\eta$, we also have
\begin{align*}
 F_\delta(x-z)
 &\geq b\bigl((\rho_K(x-z)-\delta)_+\bigr)\\
 &\geq b\bigl((\rho_K(x)-\delta-\eta)_+\bigr)\\
 &\geq b(\rho_K(x))-L(\delta+\eta).
\end{align*}
Averaging with $\eta^{-d}\varphi(\cdot/\eta)$ proves \eqref{eq:liftprofile}.
\end{proof}

\section{Evaluating the bound for a measure}\label{app:measure-integral}
The radius from Theorem~\ref{thm:smallball} is $(m/(nB))^{1/s}$.
For an explicit expression covering all parameter ranges, recall that
$b$ becomes constant at radius $R$ from \eqref{eq:radii}.
Split the integral at $\min\{(m/(nB))^{1/s},R,2L/H\}$: here the derivative changes from $Hr/2$ to
$L$, unless the interval ends first. By \eqref{eq:bderiv},
$\int_0^{(m/(nB))^{1/s}} b'(r)(m-nBr^s)\dd r$ equals
\begin{align*}
 &\frac H2\Biggl(\frac{m\min\{(m/(nB))^{1/s},R,2L/H\}^2}{2}\\
 &\qquad-\frac{nB\min\{(m/(nB))^{1/s},R,2L/H\}^{s+2}}{s+2}\Biggr)\\
 &\qquad+L\Biggl(m\left(\min\{(m/(nB))^{1/s},R\}-\min\{(m/(nB))^{1/s},R,2L/H\}\right)\\
 &\qquad-\frac{nB\left(\min\{(m/(nB))^{1/s},R\}^{s+1}-\min\{(m/(nB))^{1/s},R,2L/H\}^{s+1}\right)}{s+1}\Biggr).
\end{align*}

\section*{Data availability}
No empirical data are used in this theoretical study.

\section*{Funding}
This work received no external funding.

\section*{Declaration of competing interest}
The author declares no competing interests.

\section*{Declaration of generative AI and AI-assisted technologies in the manuscript preparation process}
Generative AI (ChatGPT, OpenAI) was used to check mathematical arguments and edit its language.

\end{document}